\documentclass[11pt,leqno,oneside,letterpaper]{amsart}

\usepackage[T1]{fontenc}
\usepackage[utf8]{inputenc}
\usepackage{lmodern}
\usepackage{eulervm}
\usepackage{amsmath,amssymb,amsthm,mathtools}
\usepackage[margin=1.15in]{geometry}
\usepackage{booktabs,array}
\usepackage{enumitem}
\usepackage{xurl}
\usepackage[colorlinks,linkcolor=blue,citecolor=blue]{hyperref}
\usepackage[nameinlink,noabbrev]{cleveref}
\usepackage{orcidlink}

\setlist{leftmargin=*,itemsep=3pt,topsep=5pt}
\newtheorem{theorem}{Theorem}[section]
\newtheorem{proposition}[theorem]{Proposition}
\newtheorem{lemma}[theorem]{Lemma}
\newtheorem{corollary}[theorem]{Corollary}
\newtheorem{conjecture}[theorem]{Conjecture}
\theoremstyle{definition}

\theoremstyle{remark}
\newtheorem{remark}[theorem]{Remark}
\numberwithin{equation}{section}
\newcommand{\Z}{\mathbb Z}
\newcommand{\Q}{\mathbb Q}
\newcommand{\OK}{\mathcal O_K}
\newcommand{\Nm}{\operatorname N}
\newcommand{\Tr}{\operatorname{Tr}}
\newcommand{\disc}{\operatorname{disc}}

\newcommand{\gp}[1]{\langle#1\rangle}

\title{Weighted unit equations in cyclic cubic orders}
\author[J.~Lu]{Junyu Lu\,\orcidlink{0009-0001-4973-3198}}
\address{Department of Mathematics, Yang-En University, Majia Town, Luojiang District, Quanzhou City, Fujian Province, China, 362014}
\email{jy.lu@outlook.com}
\date{}
\hypersetup{pdftitle={Weighted unit equations in cyclic cubic orders},pdfauthor={Junyu Lu},pdfsubject={Effective bounds and the consecutive classification}}
\begin{document}
\begin{abstract}
    We study nonrational integral units $u,v$ satisfying $au+sv=1$, with integers $1<a<s$, in varying cyclic cubic fields. The ratio $u/v$ generates the joint order $R=\Z[u,v]$ and gives effective finiteness for each fixed first weight and each fixed gap. Consecutive weights yield exactly seven simultaneous-isomorphism classes, at $a=2,4,13$; in every proper order among them, $\{-1,u,v\}$ generates $R^\times$.
\end{abstract}
\maketitle
\enlargethispage{3pt}

\section{Introduction}\label{sec:intro}
Let $K$ be a cyclic number field of degree $3$ and $\OK$ its ring of integers. For integers $1<a<s$, we consider the equation
\begin{equation}\label{eq:weighted}
    au+sv=1,
\end{equation}
where $u,v\in\OK^\times\setminus\Q$. Following Evertse--Gy\H{o}ry~\cite[\S0]{evertsegyory1988}, we call~\eqref{eq:weighted} a weighted unit equation: $a,s$ are its fixed coefficients, or weights, and the unknowns belong to the finitely generated multiplicative group $\OK^\times$. Such equations connect the multiplicative structure of unit groups with additive relations in number fields.

Setting $t=u/v$, the equation determines
\[
    R=\Z[u,v]=\Z[t],\qquad H=\gp{-1,u,v}=\gp{-1,t,at+s};
\]
the equalities and the reconstruction from $t$ are proved in Lemma~\ref{lem:ratio}. We distinguish the additive index $J=[\OK:R]$ from the unit indices
\[
    M=[\OK^\times:H],\qquad I_R=[R^\times:H].
\]
We study the order $R$ and ask whether $H=R^\times$, that is, $I_R=1$. We prove finiteness of the indices $M$ and $I_R$ using unit-generation theorems or the exponents of the units in a known basis. Note that the equality $R=\OK$ alone does not imply $H=R^\times$.

\subsection{Conventions}
The field $K$ is totally real. We write $\Nm=\Nm_{K/\Q}$ and $\Tr=\Tr_{K/\Q}$, and denote its field discriminant by $D_K=\disc(\OK)$ and its integer conductor by $\mathfrak f_K\in\Z_{>0}$. For a polynomial $f$, $\disc(f)$ denotes its polynomial discriminant. Write $n=\Nm(u)$, $T=\Tr(u)$, $e=\Nm(v)$, and $c=ne=\Nm(u/v)$, with $n,e\in\{\pm1\}$. And write the minimal polynomial of $u$ as $f_u(X)=X^3-TX^2+BX-n$.

We regard $(K,u,v)$ and $(K',u',v')$ as equivalent if there is a $\Q$-isomorphism sending $K$ to $K'$, $u$ to $u'$ and $v$ to $v'$. In a fixed cyclic field, each orbit consists of three simultaneous conjugates. We do not identify pairs by sign changes, inversion, or interchange; the condition $a<s$ fixes the ordering. When a result also allows $s<a$, we state this explicitly and keep $a$ as the first weight.

If either unit is rational, both are $\pm1$; the sole possibility under $1<a<s$ is $(-1,1)$ with $s=a+1$, excluded throughout. A cubic with positive square discriminant defines a cyclic field only after irreducibility is established. The consecutive classification uses finite integer checks after the square-gap reductions; their domains and surviving values are displayed in Section~\ref{sec:consecutive}. The unit-generation arguments below use explicit identities and cited theorems, without a computed ambient unit group.

\subsection{Main results}
Fixing the gap between the weights gives effective bounds even as the cyclic cubic field varies.

\begin{samepage}
    \begin{theorem}\label{thm:fixedgap}
        For every fixed integer $d\ge2$, the solutions of~\eqref{eq:weighted} with $s=a+d$, over all cyclic cubic fields, fall into finitely many classes under simultaneous $\Q$-isomorphism of the field and the ordered pair $(u,v)$. Representatives of these classes can be effectively determined from $d$. If $X^3-TX^2+BX-n$ is the minimal polynomial of $u$, with $n=\Nm(u)$, then
        \[
            2\le a\le d^3-d+1,\qquad s\le d^3+1,\qquad
            |B|<64(d^3+1)^7.
        \]
    \end{theorem}
\end{samepage}

Corollary~\ref{cor:fixed} also gives effective finiteness when $a$ is fixed, without assuming an upper bound on $s$. We prove Theorem~\ref{thm:fixedgap} in Section~\ref{sec:fixedgap}. The argument does not apply when the gap is one, so we treat consecutive weights separately.

\begin{theorem}\label{thm:consecutive}
    For integers $a\geq2$, suppose $K$ is a cyclic cubic number field and  $u,v\in\OK^\times\setminus\Q$ with
    \[
        au+(a+1)v=1.
    \]
    Up to simultaneous $\Q$-isomorphism of $K$ and the ordered pair $(u,v)$, the seven rows below are the complete list. In every row $\Nm(u)=-1$ and
    \[
        f_u(X)=X^3-TX^2+BX+1.
    \]
    Here $e=\Nm(v)$, $r=\Nm(u+v)$, $\mathfrak f_K$ is the field conductor, and $J=[\OK:\Z[u,v]]$.
    \begin{center}
        \begin{tabular}{rrrrrrr}
            \toprule
            $a$ & $e$  & $r$  & $T$      & $B$     & $\mathfrak f_K$ & $J$  \\
            \midrule
            2   & $-1$ & $-1$ & $-12$    & $-15$   & 7               & 1    \\
            2   & $-1$ & 1    & 24       & 3       & 9               & 1    \\
            2   & 1    & 1    & 15       & 12      & 7               & 1    \\
            2   & 1    & 9    & 159      & 84      & 9               & 57   \\
            4   & 1    & $-3$ & $-303$   & $-72$   & 9               & 17   \\
            4   & 1    & $-1$ & $-103$   & $-22$   & 7               & 1    \\
            13  & 1    & $-7$ & $-17839$ & $-1369$ & 7               & 1247 \\
            \bottomrule
        \end{tabular}
    \end{center}
    Each row means: take any root $u$ of the indicated polynomial and put $v=(1-au)/(a+1)$. The sole rational solution is $(u,v)=(-1,1)$, present for every $a\geq2$ in every cyclic cubic field.
\end{theorem}

Thus consecutive weights occur only at $a=2,4,13$, and four of the seven classes have unit sum. Section~\ref{sec:consecutive} proves Theorem~\ref{thm:consecutive}: square gaps exclude unbounded parameter ranges, a difference-of-squares reduction settles the remaining finite range, and explicit witnesses establish existence.

For fixed weights, the norm equation for $at+s$ gives a linear relation between the coefficients of the minimal polynomial of $t$. We then bound the coefficients for which its discriminant can be a square. For consecutive weights, the minimal polynomial of $u+v$ gives a simpler calculation. Proposition~\ref{prop:consecutive-generation} proves that every proper joint order in the consecutive list satisfies $R^\times=H$. We use Thomas's and Minemura's theorems for two of these orders and a residue homomorphism for the third.

The weights $(7,9)$ give an example with $R=\OK$ and $H\ne R^\times$. In the conductor-$19$ field, the units in Section~\ref{sec:order-example} satisfy $J=1$ and $I_R=M=4$. Thus generation can fail in a maximal joint order, even though it holds in all three proper consecutive orders. The orders generated by the individual units also differ from $R$: Proposition~\ref{prop:orderindices} gives $[\OK:\Z[u]]=s^3J$ and $[\OK:\Z[v]]=a^3J$.

We also study the effect of Gras's trace identity on the joint order. If $\Nm(u)u$ or $\Nm(v)v$ satisfies that identity, then $R$ is Galois invariant (Proposition~\ref{prop:gras-invariance}). The index of the order generated by that unit is $J$ times the cube of the coefficient of the other unit. For unit generation, Theorem~\ref{thm:reduction} describes the coefficients not covered by Minemura's criterion as three bounded intervals on the coefficient line. When $s=a^3\pm1$, only fourteen polynomial families remain, with $a$ still variable (Proposition~\ref{prop:extreme-reduction}). Unit generation in all proper joint orders and a uniform bound on the weights remain open; we state these conjectures in Section~\ref{sec:global}.

\subsection{Relation to earlier work}
Brindza--Gy\H{o}ry~\cite[Theorem~1 and Corollary~1]{brindzagyory} imply ineffective finiteness of the positive normalized triples $(a,s,1)$ that admit a nonrational unit solution in a fixed cyclic cubic field. The absence of a quadratic subfield makes the pairs nontrivial in their terminology, and normalization fixes the rational-equivalence ambiguity. Our fixed-weight and fixed-gap results are effective and allow the field to vary. The square-gap argument is an elementary instance of Runge's method; see Levin~\cite[Theorem~1.2]{levin}.

Vukusic--Ziegler~\cite[\S2, Theorem~1]{vukusicziegler2022} solve $\varepsilon+\delta=m$ for units of the equation orders $\Z[\rho]$, where $\rho^3-b\rho^2-(b+3)\rho-1=0$, for integers $b\ge-1$ and $|m|\le\max\{1,|b|^{1/3}\}$. Their proof uses Thomas's fundamental units in $\Z[\rho]$; we use the same theorem for members of this family in Sections~\ref{sec:order-example}, \ref{sec:consecutive-units}, and~\ref{sec:prime-unit-index}.

Tinkov\'a--Visser--Yatsyna~\cite[Theorem~1.2]{tinkovavisseryatsyna} classify ordinary unit sums $\varepsilon+\delta=m$ in cyclic cubic fields. Their equivalence permits sign changes and interchange, whereas our ordered pairs are identified only by simultaneous field isomorphism. The terms $au$ and $sv$ in~\eqref{eq:weighted} are not units. When $w=u+v$ is a unit, however, $u/w+v/w=1$ connects the two problems. The $r=\pm1$ consecutive cases also enter Komatsu's classification~\cite[Corollary~3.3 and Proposition~3.4]{komatsu}; Section~\ref{sec:thomas} records the substitutions.

Thomas--Vasquez~\cite[Theorem~2.14]{thomasvasquez} parametrize complete pairs of cubic units by twelve integral coordinates. Taking $\alpha=nu$ and $\beta=e/v$ gives a complete pair with first coordinate block $(0,n,nes,a)$. Section~\ref{sec:complete-pairs} checks the basis conditions and shows that this change preserves $R$ and $H$. We then impose the weighted equation and study the units of the resulting joint order.

For unit generation we use Thomas's Theorems~(1.4)(iii) and~(3.10)~\cite{thomas1979} and Minemura's Theorem~2~\cite{minemura}. These theorems concern equation orders. Gras's monogeneity criterion~\cite[\S I]{gras1976} motivates the trace identity in Section~\ref{sec:trace}. Our proofs require neither a classification of Galois-invariant orders nor a computation of the full unit groups of their maximal orders.

\section{The ratio order and the coefficient line}\label{sec:ratio}
Let $t=u/v$. Since $u$ and $v$ are integral units, both $t$ and $at+s=v^{-1}$ are integral units. Taking the norm of the latter identity gives $\Nm(at+s)=\Nm(v)$, a linear relation between the two nonconstant coefficients of the minimal polynomial of $t$. The following lemma shows that this relation also suffices to reconstruct $u$ and $v$ as integral units.
\begin{lemma}\label{lem:ratio}
    Every pair in \eqref{eq:weighted} satisfies
    \begin{equation}\label{eq:filters}
        \gcd(a,s)=1,\qquad s\mid na^3-1,\qquad a\mid es^3-1.
    \end{equation}
    The ratio $t=u/v$ is a nonrational integral unit, with minimal polynomial
    \begin{equation}\label{eq:P}
        P(X)=X^3+EX^2+FX-c,\qquad E,F\in\Z,\quad c=ne,
    \end{equation}
    and
    \begin{equation}\label{eq:normline}
        s^3-as^2E+a^2sF+ca^3=e.
    \end{equation}
    Moreover
    \begin{equation}\label{eq:R}
        R=\Z[t],\qquad H=\gp{-1,t,at+s}.
    \end{equation}
    Conversely, let $t$ be a root of an irreducible polynomial \eqref{eq:P} with positive square discriminant satisfying \eqref{eq:normline}. Then $v=(at+s)^{-1}$ and $u=tv$ give a pair satisfying~\eqref{eq:weighted}.
\end{lemma}
\begin{proof}
    A positive common divisor $d$ of $a,s$ would divide $1$ in $\OK$. Then $1/d$ would be a rational algebraic integer, forcing $d=1$. Taking norms of $au=1-sv$ modulo $s$, and of $sv=1-au$ modulo $a$, proves the two divisibilities. Inverses of integral units are integer polynomials in those units, so $t$ is integral and invertible. If $t$ were rational, $v(at+s)=1$ would force both units to be rational. Thus $[\Q(t):\Q]=3$. The norm of $t$ is $c$, which gives \eqref{eq:P}. Taking norms of $at+s=v^{-1}$ gives \eqref{eq:normline}.

    We have $t=uv^{-1}\in R$. Conversely, $at+s\in\Z[t]$ has norm $e=\pm1$, so its inverse $v$ belongs to $\Z[t]$, as does $u=tv$. Hence $R=\Z[t]$. The identities $t=uv^{-1}$ and $at+s=v^{-1}$ express each pair of generators in terms of the other; adjoining $-1$ gives the second equality in \eqref{eq:R}.

    For the converse, a monic irreducible cubic of positive square discriminant defines a totally real cyclic cubic field. Its constant coefficient makes $t$ a unit, and \eqref{eq:normline} makes $at+s$ a unit. Thus $v=(at+s)^{-1}$ and $u=tv$ are integral units of norms $e$ and $ce=n$, respectively, and satisfy the weighted equation. Starting from a pair and reconstructing it from its ratio returns the same pair. Starting from a root $t$ returns a pair with ratio $t$. This proves the correspondence of simultaneous-isomorphism classes.
\end{proof}

For fixed weights $a,s$ and norm signs $n,e$, the norm identity leaves one integer parameter in the coefficients of $P$. Put
\begin{equation}\label{eq:t}
    \ell=\frac{s^3-e+ca^3}{as}.
\end{equation}
The two norm congruences make the numerator of $\ell$ divisible by both $a$ and $s$. Since $\gcd(a,s)=1$, we have $\ell\in\Z$. Equation~\eqref{eq:normline} is equivalent to $sE-aF=\ell$. Choose the unique $E_0\in\{0,\ldots,a-1\}$ with $sE_0\equiv\ell\pmod a$, and set $F_0=(sE_0-\ell)/a$. Every integral solution is then of the form
\begin{equation}\label{eq:line}
    E=E_0+ak,\qquad F=F_0+sk,\qquad k\in\Z.
\end{equation}
Every integer $k$ gives integral coefficients satisfying the norm equation. Substituting these coefficients into the discriminant gives the quartic
\begin{equation}\label{eq:disc}
    D(k)=E^2F^2-4F^3+4cE^3-18cEF-27,
\end{equation}
whose leading coefficient is $(as)^2$. By Lemma~\ref{lem:ratio}, we must find the positive square values of $D(k)$ for which $P$ is irreducible. Section~\ref{sec:lattice} bounds the integers $k$ that need to be considered.

To express the reduction in the coefficients of $u$, write its minimal polynomial as $X^3-TX^2+BX-n$ and define $\gamma=(1-a^3n-s^3e)/a\in\Z$. Then
\begin{equation}\label{eq:old}
    \begin{aligned}
        T & =aB+\gamma,\qquad F=\frac{B-3an}{es},       \\
        B & =es^2k+esF_0+3an,\qquad \disc(f_u)=s^6D(k).
    \end{aligned}
\end{equation}

The ratio gives a generator for the joint order $R$. To compare $R$ with the orders generated by the individual units, we use the following index formulas.
\begin{samepage}
    \begin{proposition}\label{prop:orderindices}
        For every pair~\eqref{eq:weighted},
        \[
            [R:\Z[u]]=s^3,\qquad [R:\Z[v]]=a^3.
        \]
        Consequently, if $J=[\OK:R]$, then
        \[
            [\OK:\Z[u]]=s^3J,\quad [\OK:\Z[v]]=a^3J,
            \quad J=\gcd([\OK:\Z[u]],[\OK:\Z[v]]).
        \]
    \end{proposition}
\end{samepage}
\begin{proof}
    Choose integers $u_0,v_0$ with $au_0+sv_0=1$, and let $z=v_0u-u_0v$. Substitution gives
    \[
        u=u_0+sz,\qquad v=v_0-az,
    \]
    so $R=\Z[z]$. The determinants of the coordinate changes from $(1,z,z^2)$ to $(1,u,u^2)$ and $(1,v,v^2)$ are $s^3$ and $-a^3$. Since $z$ has degree three, these are full integral bases of their respective equation orders, proving the first identities. Multiplying by $J$ and using $\gcd(a,s)=1$ proves the rest.
\end{proof}

\subsection{A maximal joint order with unit index four}\label{sec:order-example}
The additive index $J$ does not determine the unit index $I_R$. To see this, let $t$ be a root of
\[
    P(X)=X^3-2X^2-5X-1.
\]
The values $P(1)=-7$ and $P(-1)=1$ prove irreducibility, and $\disc(P)=19^2$ proves that $K=\Q(t)$ is cyclic. If $\Z[t]$ were not maximal, its index would be $19$ and the field discriminant would be $1$, contradicting Minkowski's bound $D_K\ge81/4$ for a totally real cubic field. Hence $\Z[t]=\OK$. Thomas's Theorem~(3.10)~\cite{thomas1979}, with parameter $2$, gives
\[
    \OK^\times=\gp{-1,t,1+t}.
\]
The defining relation $t^3=2t^2+5t+1$ gives
\[
    (1+t)^4=23t^2+35t+7=t^2(7t+9).
\]
Thus
\[
    u=t^3(1+t)^{-4},\qquad v=t^2(1+t)^{-4}
\]
are integral units satisfying $7u+9v=1$, and $u/v=t$. Their norms are one because $\Nm(t)=1$ and $\Nm(1+t)=-1$. In the basis $t,1+t$ of $\OK^\times/\{\pm1\}$, their exponent columns form
\[
    \begin{pmatrix}3&2\\-4&-4\end{pmatrix},\qquad
    \left|\det\begin{pmatrix}3&2\\-4&-4\end{pmatrix}\right|=4.
\]
Consequently $J=1$ and $M=I_R=4$: the joint order is maximal, but $H$ has index four in its unit group. Proposition~\ref{prop:orderindices} also gives $[\OK:\Z[u]]=9^3$ and $[\OK:\Z[v]]=7^3$. These additive indices describe the individual equation orders and are distinct from the unit index.

\subsection{The complete-pair coordinates}\label{sec:complete-pairs}
To place the ratio description in the coordinate construction of Thomas--Vasquez~\cite[\S2]{thomasvasquez}, we first check their basis conditions. They call a pair $(\alpha,\beta)$ complete when both norms are one and both $\{1,\alpha,\beta\}$ and $\{\alpha\beta,\alpha,\beta\}$ are $\Q$-bases. For our pair set
\[
    \alpha=nu,\qquad \beta=e/v.
\]
Both norms are one, and $\beta=es/(1-an\alpha)$. If $A+B\alpha+C\beta=0$ with rational coefficients, multiplication by $1-an\alpha$ yields
\[
    -anB\alpha^2+(B-anA)\alpha+A+Ces=0.
\]
Since $\alpha$ has degree three, the coefficients successively give $B=A=C=0$. The identity
\[
    a\alpha\beta=n\beta-nes
\]
then proves the second basis condition, because its constant term is nonzero. In the notation of their equation~(2.3), this is the normalized first block $(x_1,x_2,x_3,x_4)=(0,n,nes,a)$: its entries have greatest common divisor one and $x_4>0$.

We recover the original pair by $u=n\alpha$ and $v=e\beta^{-1}$. Both changes use integer polynomials in units and their inverses, so they preserve $R$ and $H$. Thus our pairs belong to the part of the Thomas--Vasquez parametrization with first block $(0,n,nes,a)$. For fixed weights, the coefficient line~\eqref{eq:line} gives a simpler description in terms of $E$ and $F$. We still need to determine when the discriminant is a positive square and whether $H=R^\times$.

\section{Bounding square values of the discriminant}\label{sec:lattice}
The leading coefficient of $D(k)$ is a square. We can therefore choose a quadratic polynomial $Q(k)$ such that $D(k)-Q(k)^2$ is linear. Suppose there is a constant $\delta>0$ such that $|Q(k)-y|\ge\delta$ for all integers $k,y$. For large $|k|$, the distance from $Q(k)^2$ to an integer square then grows quadratically, while $D(k)-Q(k)^2$ grows at most linearly. This gives a bound on $|k|$ whenever $D(k)$ is a square. We first prove the bound and then show that the required separation holds for the weighted equation.
\begin{lemma}\label{lem:gap}
    Let $D(k)\in\Z[k]$ be a quartic with leading coefficient $A^2$, $A\in\Z_{>0}$. Match its three highest coefficients to write uniquely
    \begin{equation}\label{eq:Q}
        Q(k)=Ak^2+Lk+C_Q,\qquad D(k)-Q(k)^2=hk+j,
    \end{equation}
    with $L,C_Q,h,j\in\Q$. For $b\in\Q$, let $\operatorname{den}(b)$ denote its positive denominator in lowest terms. Define
    \[
        \begin{gathered}
            q=\operatorname{lcm}(\operatorname{den}L,\operatorname{den}C_Q),\qquad G=\gcd(q,qL),\\
            r_0\equiv qC_Q\pmod G,\qquad 0\le r_0<G.
        \end{gathered}
    \]
    If $d_0=\min(r_0,G-r_0)>0$, set $\delta=d_0/q$. For all integers $k,y$, $|Q(k)-y|\ge\delta$. Consequently $D(k)$ cannot be a square when
    \begin{equation}\label{eq:exclude}
        Q(k)\ge\delta,\qquad 2\delta Q(k)-\delta^2>|hk+j|.
    \end{equation}
\end{lemma}
\begin{proof}
    Since $A$ is an integer, $qQ(k)-qy\equiv qC_Q\pmod G$. Every integer in this congruence class has absolute value at least $d_0$. Dividing by $q$ gives the claimed lower bound for $|Q(k)-y|$.

    Suppose $y\ge0$ is an integer and $Q\ge\delta$. If $y\le Q$, let $z=Q-y\in[\delta,Q]$; the function $z(2Q-z)$ is increasing on this interval. Thus $|Q^2-y^2|\ge\delta(2Q-\delta)$. If $y\ge Q$, the same bound follows from $(y-Q)(y+Q)\ge\delta(2Q+\delta)$. The strict residual inequality in \eqref{eq:exclude} contradicts $D(k)=y^2$.
\end{proof}

The bound can be computed using only rational arithmetic. Set
\begin{equation}\label{eq:HK}
    H_0=|L|+\frac{|h|}{2\delta},\qquad
    K_0=|C_Q|+\frac{|j|}{2\delta}+\frac\delta2.
\end{equation}
Choose an integer $N>0$ such that both quadratic expressions below and their derivatives are positive at $z=N$:
\begin{equation}\label{eq:bounds}
    Az^2-|L|z-|C_Q|-\delta>0,\qquad Az^2-H_0z-K_0>0.
\end{equation}
Such an $N$ exists because both leading coefficients are positive. The expressions then remain positive for every real $z\ge N$. Since $Q(k)\ge A|k|^2-|L||k|-|C_Q|$, \eqref{eq:exclude} holds for all $|k|\ge N$. Every square value therefore occurs in the finite interval $|k|<N$.

We now verify the required separation $d_0>0$ for the quartics arising from the weighted equation.
\begin{proposition}\label{prop:positive}
    For every choice of weights and norm signs satisfying \eqref{eq:filters} with $a,s>1$, the quartic \eqref{eq:disc} on the coefficient line \eqref{eq:line} has $d_0>0$ in Lemma~\ref{lem:gap}. Thus the lemma gives a finite bound for these families, even when the quartic has repeated roots.
\end{proposition}
\begin{proof}
    In terms of the coefficient $B$ in~\eqref{eq:old}, the quadratic approximation is
    \begin{align}\label{eq:QB}
        Q_B(B)={} & aB^2+(\gamma-2/a-2na^2)B\notag        \\
                  & -4na\gamma+2\gamma/a^2+5n-2/a^3-2a^3.
    \end{align}
    Direct coefficient matching gives a linear residual. For integral $B,\gamma$,
    \[
        a^3 Q_B(B)\equiv-2\pmod a.
    \]
    If $a>2$, this prevents $Q_B(B)$ from being integral. For $a=2$ its nonconstant coefficients are integral, while its constant coefficient is
    \[
        -\frac{15(3+2\gamma)}4\quad(n=1),\qquad
        -\frac{17(5-2\gamma)}4\quad(n=-1),
    \]
    neither of which is integral. The discriminant identity in~\eqref{eq:old} and the positive leading coefficients of the two quadratic approximations give
    \begin{equation}\label{eq:scale}
        Q_B(es^2k+esF_0+3an)=s^3Q(k).
    \end{equation}
    Thus $Q(k)$ is never integral at an integral $k$.

    As the quadratic coefficient $A=as$ is integral, $Q(k)$ is integral for some integer $k$ if and only if the linear congruence $qLk+qC_Q\equiv0\pmod q$ is solvable. This occurs precisely when $G\mid qC_Q$, equivalently $d_0=0$. The nonintegrality just proved gives $d_0>0$, without any assumption that the quartic is squarefree.
\end{proof}

The residual $hk+j$ has no integer root. Indeed, if $hk+j=0$ for an integer $k$, then $D(k)=Q(k)^2$ would be an integer, forcing the rational number $Q(k)$ to be an integer. In the $B$ coordinate, the linear coefficient of the residual factors as
\[
    h_B=\frac{4(1+na^3)(1+es^3)(na^3-es^3)}{a^4}\ne0.
\]
Nonvanishing follows from $a,s>1$ and coprimality.

\begin{remark}
    For an arbitrary quartic, $d_0$ may be zero, so the bound above may not apply. When $Q\ne y$, one can still use $|Q-y|\ge1/q$, but integer roots of $hk+j$ must be considered separately. If $h=j=0$, there may be infinitely many square values, as in $D(k)=k^4$. Proposition~\ref{prop:positive} rules out these difficulties for the weighted families considered here.
\end{remark}

For a fixed first weight, the norm congruences leave only finitely many second weights. Applying the discriminant bound to each of them gives the following classification procedure.
\begin{corollary}\label{cor:fixed}
    Fix an integer $a\ge2$. The nonrational integral-unit solutions of $au+sv=1$, with integer $s>1$, over all cyclic cubic fields can be classified by a finite exact calculation. Here $a$ remains the first weight, and $s<a$ is allowed.
\end{corollary}
\begin{proof}
    For each norm sign $n$, the weight $s$ is a positive divisor greater than one of the nonzero integer $|na^3-1|$. Retain the norm signs satisfying \eqref{eq:filters}. For each resulting choice, \eqref{eq:line} and Proposition~\ref{prop:positive} give a finite interval of $k$ containing every square discriminant. Check the integers in this interval, retain positive square discriminants, and discard reducible cubics. For a monic cubic with constant coefficient $\pm1$, the only possible rational roots are $1$ and $-1$. Lemma~\ref{lem:ratio} reconstructs the pairs and proves that none are omitted. Distinct values of $k$ for fixed weights and signs give distinct minimal polynomials of the ratio, hence distinct orbits of ordered pairs.
\end{proof}
The corollary gives a finite procedure for each fixed $a$. It gives no bound on the running time that would make the procedure practical for all $a$, and the discriminant calculation alone does not prove unit generation.

\section{An effective bound for every fixed gap}\label{sec:fixedgap}
\begin{proof}[Proof of Theorem~\ref{thm:fixedgap}]
    We first bound the weights, then the coefficient $B$ of the minimal polynomial of $u$. Substituting $s=a+d$ in the norm congruences~\eqref{eq:filters} gives
    \begin{equation}\label{eq:gapfilters}
        s\mid d^3n+1,\qquad a\mid d^3e-1,\qquad \gcd(a,d)=1.
    \end{equation}
    Since $d\ge2$, the first dividend is nonzero. This proves the weight bounds and, together with Corollary~\ref{cor:fixed}, already proves effective finiteness. To obtain the stated uniform coefficient bound, use $\gamma=(1-a^3n-s^3e)/a$, the matched polynomial $Q_B=aB^2+LB+C_Q$ of~\eqref{eq:QB}, and its residual $h_BB+j_B$. Thus
    \[
        L=\gamma-2/a-2na^2,\quad
        C_Q=-4na\gamma+2\gamma/a^2+5n-2/a^3-2a^3,
    \]
    \[
        h_B=\frac{4(1+na^3)(1+es^3)(na^3-es^3)}{a^4},
        \qquad j_B=-4n\gamma^3-27-C_Q^2.
    \]
    Set $M_0=d^3+1$ and $q=a^3$. Elementary estimates, using $a,s\le M_0$, give
    \[
        |\gamma|\le\frac{3M_0^3}{a},\quad |L|\le\frac{6M_0^3}{a},\quad
        |C_Q|\le16M_0^3,\quad |h_B|\le\frac{32M_0^6}{a},\quad |j_B|\le47M_0^9.
    \]
    For the bound on $h_B$, the factorization above gives
    \[
        \begin{aligned}
            |h_B| & \le\frac{4(1+a^3)(1+s^3)(a^3+s^3)}{a^4} \\
                  & \le\frac{4(2a^3)(2M_0^3)(2M_0^3)}{a^4}
            =\frac{32M_0^6}{a}.
        \end{aligned}
    \]
    For the bound on $j_B$, use
    \[
        |C_Q|\le\tfrac{59}{4}M_0^3+\tfrac{21}{4}\le16M_0^3,
        \qquad |j_B|\le\tfrac{27}{2}M_0^9+256M_0^6+27\le47M_0^9.
    \]
    Let $H_0=|L|+q|h_B|$, $K_0=|C_Q|+q|j_B|$. Then
    \[
        H_0/a\le33M_0^7,\qquad K_0/a\le8M_0^3+47M_0^{11}.
    \]
    At $z=64M_0^7$,
    \[
        z^2-(H_0/a)z-K_0/a
        \ge1984M_0^{14}-47M_0^{11}-8M_0^3>0,
    \]
    and its derivative is at least $95M_0^7>0$. For all $|B|\ge64M_0^7$ we therefore have $Q_B(B)>q|h_BB+j_B|$ and $Q_B(B)>0$. If the discriminant were $y^2$ with $y\in\Z_{\ge0}$, the nonintegrality argument in Proposition~\ref{prop:positive} would give $|Q_B(B)-y|\ge1/q$. Hence
    \[
        |h_BB+j_B|=|y-Q_B(B)|\bigl(y+Q_B(B)\bigr)
        \ge\frac{Q_B(B)}q>|h_BB+j_B|,
    \]
    a contradiction. This proves the coefficient bound.

    To reconstruct the pairs from $B$, impose the congruence and trace relation
    \[
        B\equiv(3-2a\gamma)(a^2)^{-1}\pmod{s^2},\qquad T=aB+\gamma.
    \]
    Since $s\mid\gamma$, the congruence makes both $(3-aT)/s$ and $(3-2aT+a^2B)/s^2$ integral. These are the trace and quadratic coefficient of $v=(1-au)/s$. Its norm is $e$, so $v$ is an integral unit. Conversely, both coefficients are integral for every solution. It therefore remains to require positive square discriminant and to exclude the possible rational roots $\pm1$.
\end{proof}
The failure of the first step for $d=1,n=-1$ explains the need for the separate consecutive theorem. A fixed-field finiteness theorem alone would not justify Theorem~\ref{thm:fixedgap}, since the field here varies.

\section{The complete consecutive classification}\label{sec:consecutive}
Throughout this section, $s=a+1$. To prove Theorem~\ref{thm:consecutive}, we use $w=u+v$, whose minimal polynomial has a simpler discriminant than those of $u$ and $v$. We first determine the possible norm signs, then bound the parameters in the polynomial of $w$. The remaining finite checks give seven candidates, for which we construct explicit pairs.

\subsection{Norm signs and the cubic model}\label{sec:consec-signs}
With $s=a+1$, the norm congruences~\eqref{eq:filters} become
\[
    a+1\mid a^3n-1,\qquad a\mid(a+1)^3e-1.
\]
The first condition implies $a+1\mid n+1$, so $n=-1$. The second says $a\mid e-1$. Thus $e=1$ for $a\geq3$; only $a=2$ permits $e=-1$. If either unit is rational, both are in $\{\pm1\}$, and direct substitution leaves only $(-1,1)$.

Put $w=u+v$. Then
\begin{equation}\label{eq:consec-affine}
    u=(a+1)w-1,\qquad v=1-aw.
\end{equation}
For a nonrational pair, $w$ is a nonrational algebraic integer of degree three. Write its minimal polynomial as $X^3-T_wX^2+bX-r$, where $r=\Nm(w)\ne0$.

First suppose $e=1$. The two norm identities give
\[
    T_w=(a+1)b-(a+1)^2r,\qquad T_w=ab-a^2r.
\]
Consequently, with
\[
    P=a(a+1),\qquad S=2a+1,\qquad S^2=4P+1,
\]
we obtain
\begin{equation}\label{eq:consec-model}
    g_{a,r}(X)=X^3-PrX^2+SrX-r.
\end{equation}
Conversely, an irreducible polynomial \eqref{eq:consec-model} with positive square discriminant defines a cyclic cubic field; \eqref{eq:consec-affine} then gives integral units of norms $-1,1$ satisfying the equation. Thus this reduction is sufficient as well as necessary after checking irreducibility.

A direct discriminant calculation yields
\begin{equation}\label{eq:consec-disc}
    \disc(g_{a,r})=r^2D_a(r),\qquad
    D_a(r)=P^2r^2+2S(P-2)r-27.
\end{equation}
As $r\ne0$, cyclicity requires $D_a(r)$ to be a positive integer square. We first find the pairs $(a,r)$ satisfying this square condition and then check irreducibility.

\subsection{Square gaps for unbounded parameters}\label{sec:consec-gaps}
We continue with $e=1$. Comparing $D_a(r)$ with adjacent integer squares bounds $a$, apart from two explicit pairs. The comparison depends on the sign and size of $r$.
\subsubsection{The four small nonzero values of the norm}
For $r=\pm1,\pm2$, we have
\begin{align*}
    D_a(1)  & =(P+S-2)^2-32,       \\
    D_a(-1) & =(P-S-2)^2-32,       \\
    D_a(2)  & =(2P+S-1)^2-(6S+29), \\
    D_a(-2) & =(2P-S-1)^2+(6S-29).
\end{align*}
For $r=1$ and $a\geq3$, put $\mu=P+S-2\geq17$, so $0<32<2\mu-1$. For $r=-1$ the same argument applies to $\mu=P-S-2$ when $a\geq5$. Both discriminants lie strictly between $(\mu-1)^2$ and $\mu^2$.

For $r=2$, $\mu=2P+S-1$ satisfies
\[
    2\mu-1-(6S+29)=S^2-4S-33>0\qquad(a\geq4),
\]
again excluding a square. For $r=-2$ and every $a\geq2$, $\mu=2P-S-1\geq6$ and
\[
    6S-29>0,\qquad 2\mu+1-(6S-29)=(S-4)^2+11>0.
\]
Thus $D_a(-2)$ lies strictly between $\mu^2$ and $(\mu+1)^2$.

\subsubsection{Positive norm at least three}
Suppose $r\geq3$ and $a\geq11$, hence $S\geq23$. Set $\mu=Pr+S$ and $\Delta=4Sr+S^2+27$. Then $D_a(r)=\mu^2-\Delta$ and $\Delta>0$. To compare with the preceding square, observe that
\begin{align*}
    2\mu-1-\Delta
     & =(2P-4S)r-S^2+2S-28          \\
     & \geq \frac{S^2-20S-59}{2}>0.
\end{align*}
Here $2P-4S=(S^2-8S-1)/2>0$, and the last numerator is $10$ at $S=23$ and increases thereafter. Hence $(\mu-1)^2<D_a(r)<\mu^2$, so no square occurs.

\subsubsection{Negative norm at most minus three}
Write $r=-k$, where $k\geq3$, and suppose $a\geq4$, so $S\geq9$. Put
\[
    \mu=Pk-S>0,\qquad \delta=4Sk-S^2-27.
\]
Now $D_a(r)=\mu^2+\delta$, and both adjacent-square gaps are positive:
\begin{align*}
    D_a(r)-(\mu-1)^2
     & =(2P+4S)k-S^2-2S-28         \\
     & \geq\frac{S^2+20S-59}{2}>0, \\
    (\mu+1)^2-D_a(r)
     & =(2P-4S)k+S^2-2S+28>0.
\end{align*}
The last inequality uses $2P-4S\geq4$. Therefore the only possible square is $\mu^2$, forcing
\[
    \delta=0,\qquad k=\frac{S^2+27}{4S}.
\]
For $k$ to be an integer, we must have $S\mid27$. Since $S$ is odd and $S\geq9$, the only possibilities are
\begin{equation}\label{eq:consec-negative-candidates}
    (a,r)=(4,-3),\qquad (a,r)=(13,-7).
\end{equation}
Combining these comparisons, the only possible pair with $a\geq11$ is $(13,-7)$. For negative $r$, the comparison uses the positive number $\mu=P|r|-S$; no sign assumption on $Pr+S$ is needed.

\subsection{The remaining finite checks}\label{sec:consec-finite}
\subsubsection{The case \texorpdfstring{$(n,e)=(-1,1)$}{(n,e)=(-1,1)}}
It remains to check $a=2,\ldots,10$. Although $r$ is still unbounded, a difference of squares reduces each value of $a$ to finitely many factor pairs. Put
\[
    C_1=2S(P-2),\qquad C=16(P+1)^3.
\]
Equation $D_a(r)=Y^2$ is equivalent to
\[
    (2P^2r+C_1)^2-(2PY)^2=C.
\]
We therefore examine all signed factor pairs $\alpha\beta=C$ with $\beta\geq\alpha$ and retain those satisfying
\begin{equation}\label{eq:consec-factors}
    r=\frac{\alpha+\beta-2C_1}{4P^2}\in\Z,
    \qquad Y=\frac{\beta-\alpha}{4P}\in\Z_{>0}.
\end{equation}
\begin{samepage}
    Including both positive and negative factor pairs gives $420$ pairs to check. The pairs satisfying these conditions are
    \begin{center}
        \begin{tabular}{rrrrr}
            \toprule
            $a$ & $r$  & $Y$ & $\alpha$ & $\beta$ \\
            \midrule
            2   & 1    & 7   & 28       & 196     \\
            2   & 9    & 57  & 4        & 1372    \\
            4   & $-3$ & 51  & $-4116$  & $-36$   \\
            4   & $-1$ & 7   & $-756$   & $-196$  \\
            \bottomrule
        \end{tabular}
    \end{center}
\end{samepage}
No pair satisfies the conditions for $a=3,5,6,7,8,9,10$. For $a=2,\ldots,10$, respectively, the numbers $P+1$ are
\[
    7,\ 13,\ 21,\ 31,\ 43,\ 57,\ 73,\ 91,\ 111.
\]
Five of these numbers are prime and four are products of two distinct primes, so $C=16(P+1)^3$ has $20$ or $80$ positive divisors, respectively. Since $C$ is nonsquare, including both signs gives as many unordered factor pairs as positive divisors, hence $5\cdot20+4\cdot80=420$ pairs in total. Including $(13,-7)$ from \eqref{eq:consec-negative-candidates} gives five candidates in the norm $(-1,1)$ branch.

\subsubsection{The case \texorpdfstring{$(n,e)=(-1,-1)$}{(n,e)=(-1,-1)}}
Only $a=2$ is possible. The same substitution $u=3w-1$, $v=1-2w$ now gives
\[
    g_r(X)=X^3-(6r+3)X^2+(5r+1)X-r,
\]
with discriminant
\[
    D_-(r)=36r^4+4r^3+24r^2+12r+5.
\]
To bound $r$, we again approximate the discriminant by the square of a quadratic. Set
\[
    Q(r)=6r^2+\frac r3+\frac{215}{108}.
\]
Then
\begin{equation}\label{eq:consec-residual}
    D_-(r)-Q(r)^2=\frac{1729}{162}r+\frac{12095}{11664}.
\end{equation}
The residual root is $-12095/124488$, which is not integral. For $x=|r|$,
\[
    Q(r)-108|D_-(r)-Q(r)^2|
    \geq 6x^2-1153x-\frac{6155}{54}.
\]
At $x=193$, the right side is $45955/54>0$ and its derivative is $1163>0$, so it remains positive for $x\geq193$. On the other hand, if $D_-(r)=y^2$ with $y\in\Z_{\geq0}$, then $108Q(r)\in\Z$ and the nonzero residual implies
\[
    |Q(r)-y|\geq\frac1{108},\qquad
    |D_-(r)-Q(r)^2|\geq\frac{Q(r)}{108},
\]
a contradiction. Every square case therefore satisfies $|r|<193$.

Checking the $385$ integers $-192\leq r\leq192$ gives only $r=-1$ and $r=1$, with discriminants $49$ and $81$. Their polynomials are $X^3+3X^2-4X+1$ and $X^3-9X^2+6X-1$.

\subsection{Existence and additive order indices}
For each surviving polynomial of $w$, the rational-root test over all $\pm d$ with $d\mid|r|$ proves irreducibility. Testing only $\pm1$ would not suffice when $|r|>1$. All seven discriminants are positive squares, hence the fields are cyclic cubic.

We construct all seven classes in the fields defined by
\[
    g_7(X)=X^3-X^2-2X+1,\qquad g_9(X)=X^3-3X-1.
\]
In each row below, let $\theta$ be a root of $g_{\mathfrak f_K}$ and take $w=w_0+w_1\theta+w_2\theta^2$. The formulas~\eqref{eq:consec-affine} then give $u$ and $v$.
\begin{center}
    \begin{tabular}{rrrrrrr}
        \toprule
        $a$ & $e$  & $r$  & $\mathfrak f_K$ & $w_0$ & $w_1$  & $w_2$  \\
        \midrule
        2   & $-1$ & $-1$ & 7               & 1     & $-1$   & $-1$   \\
        2   & $-1$ & 1    & 9               & 9     & 1      & $-3$   \\
        2   & 1    & 1    & 7               & 5     & 1      & $-2$   \\
        2   & 1    & 9    & 9               & 58    & 7      & $-20$  \\
        4   & 1    & $-3$ & 9               & 10    & 23     & $-15$  \\
        4   & 1    & $-1$ & 7               & 3     & $-4$   & $-5$   \\
        13  & 1    & $-7$ & 7               & 171   & $-247$ & $-308$ \\
        \bottomrule
    \end{tabular}
\end{center}
Polynomial reduction modulo $g_{\mathfrak f_K}$ verifies each defining cubic and the identities \eqref{eq:consec-affine}, including the unit norms. To identify the field conductors and compute $J$, we first check that $1,\theta,\theta^2$ is an integral basis. The polynomials $g_7$ and $g_9$ have discriminants $49$ and $81$; their translates
\[
    g_7(X-2)=X^3-7X^2+14X-7,\qquad
    g_9(X+1)=X^3+3X^2-3
\]
are Eisenstein at $7$ and $3$. At the respective prime, the translated root is a uniformizer in a totally ramified cubic extension, and its powers form a local integral basis. Translation preserves the generated order, so the original power bases are locally maximal at the only possible index primes. Thus $49$ and $81$ are the field discriminants, and the cyclic conductor--discriminant formula gives conductors $7$ and $9$.

Since $w=u+v$ and \eqref{eq:consec-affine} holds, $R=\Z[u,v]=\Z[w]$. Basis changes yield
\[
    [R:\Z[u]]=(a+1)^3,\quad [R:\Z[v]]=a^3,
\]
and hence
\[
    [\OK:\Z[u]]=(a+1)^3J,\qquad
    [\OK:\Z[v]]=a^3J.
\]
The determinant of the coordinates of $1,w,w^2$ in the integral basis $1,\theta,\theta^2$ gives the values $J$ in Theorem~\ref{thm:consecutive}. These are additive indices of orders.

\begin{proof}[Proof of Theorem~\ref{thm:consecutive}]
    Section~\ref{sec:consec-signs} restricts the norm signs to $n=-1$, with $e=1$ unless $a=2$. The square gaps of Section~\ref{sec:consec-gaps} and the exhaustive checks of Section~\ref{sec:consec-finite} leave precisely the seven displayed classes. The witnesses above establish existence and identify their conductors and additive indices. Since $u=(a+1)w-1$ and $v=1-aw$, each defining polynomial of $w$ determines one simultaneous-isomorphism class of ordered pairs. The seven rows have distinct defining data, and the rational case was settled in Section~\ref{sec:consec-signs}.
\end{proof}

Each row of Theorem~\ref{thm:consecutive} is one Galois orbit of three ordered pairs in a fixed model of its field. Summing over all $a\geq2$, the conductor-$7$ and conductor-$9$ fields therefore have twelve and nine nonrational ordered pairs, respectively. No quotient by sign or inversion is taken.

\subsection{Relation with the ratio model and the Thomas family}\label{sec:thomas}
To study unit generation for the classified pairs, we return to the ratio $t=u/v$. When the norm signs are $(-1,1)$, its minimal polynomial is
\[
    P_t(X)=(X+1)^3+rX(aX+a+1).
\]
Indeed, $w=(t+1)/(at+a+1)$, and substitution in the minimal polynomial of $w$ gives the displayed polynomial. This identifies the consecutive calculation with the ratio equation of Section~\ref{sec:ratio}.

For negative $r=-k$, put $\lambda=-w$. Then
\[
    \phi(X)=X^3-k(aX+1)((a+1)X+1)
    =X^3-(R_0X+1)(pX+k)
\]
with $R_0=a+1$, $p=ak$. Thomas~\cite[Theorem~(1.4)(iii)]{thomas1979} assumes $R_0\ge4$, $k\ge2$, $|p|\ge2$, and either
\[
    p\le\min(R_0k-5,(R_0-1)k)
    \quad\hbox{or}\quad
    p\ge\max(R_0k+5,(R_0+1)k).
\]
When $p=mk$ and $m\le R_0-1$, Thomas's theorem gives the fundamental pair $\{m\lambda+1,R_0\lambda+1\}$ in $\Z[\lambda]$. For $a\ge3$ and $k\ge5$, the first inequalities hold with $m=a$. In particular, $(a,r)=(13,-7)$ gives $(R_0,p,k,m)=(14,91,7,13)$ and the pair $\{v,-u\}$. Since $\Z[\lambda]=R$, the theorem proves $R^\times=H$ for this row. Its conclusion concerns the units of $R$ and does not assert generation in $\OK$.

In the norm-$(-1,1)$ branch, the $r=\pm1$ subcases also lie in the unit-sum reductions associated with Komatsu~\cite[Corollary~3.3 and Proposition~3.4]{komatsu}: the inverse $1/w$ has polynomial $X(X-a)(X-a-1)-r$. In his notation $g_{s_0,t_0}(Y)=Y(Y+s_0)(Y+t_0)-1$, take $Y=1/w$ and $(s_0,t_0)=(-a,-a-1)$ when $r=1$, and $Y=-1/w$ and $(s_0,t_0)=(a,a+1)$ when $r=-1$. The cited results give $a=2$ and $a=4$, respectively, both in the conductor-$7$ field. The values $r=9,-3,-7$ lie outside that reduction. In particular, the $a=13$ row has $\Nm(u+v)=-7$.

\subsection{Unit generation in the proper orders}\label{sec:consecutive-units}
Theorem~\ref{thm:consecutive} contains three proper joint orders, with $(a,r)=(2,9),(4,-3),(13,-7)$. We prove $I_R=1$ for each of them. The first follows from Minemura's criterion in the form stated in Section~\ref{sec:minemura}, and the third from Thomas's theorem as applied above. For the middle case, we describe $R$ modulo $17$ and compare the indices of $H$ and $R^\times$ in $\OK^\times$.

\begin{proposition}\label{prop:consecutive-generation}
    Every nonrational consecutive pair with $R\subsetneq\OK$ satisfies $R^\times=\gp{-1,u,v}$.
\end{proposition}
\begin{proof}
    For $(a,r)=(2,9)$, the ratio polynomial in Section~\ref{sec:thomas} is $X^3+21X^2+30X+1$. The conditions in~\eqref{eq:mine} become
    \[
        30\pm21\ne1,\qquad \frac{105}{4}>\frac{43}{2},
        \qquad \frac{39}{2}>8,\qquad 18>\frac52.
    \]
    Minemura's theorem therefore gives generation. For $(a,r)=(13,-7)$, the substitution into Thomas's Theorem~(1.4)(iii) already established in Section~\ref{sec:thomas} gives the same conclusion.

    It remains to consider $(a,r)=(4,-3)$. We first identify $R$ inside $\OK$. Let $\theta$ be a root of $\theta^3-3\theta-1=0$. The preceding integral-basis calculation gives $\OK=\Z[\theta]$, with conductor $9$. Thomas's Theorem~(3.10) gives $\OK^\times=\gp{-1,\theta,\theta+1}$. Put
    \[
        u=\theta^6(\theta+1)^{-5},\qquad v=\theta^4(\theta+1)^{-6},
        \qquad t=u/v=\theta^2+3\theta+1.
    \]
    Reduction using $\theta^3=3\theta+1$ gives
    \[
        (\theta+1)^6=\theta^4\bigl(4\theta^2(\theta+1)+5\bigr),\qquad
        t^3-9t^2-12t+1=0.
    \]
    Hence $4u+5v=1$, and $t$ has the ratio polynomial of the required row. In the basis $1,\theta,\theta^2$,
    \[
        t^2=7+25\theta+14\theta^2,\qquad
        [\OK:\Z[t]]=\left|\det\begin{pmatrix}3&25\\1&14\end{pmatrix}\right|=17.
    \]
    Both $t$ and $t^2$ satisfy the congruence $B-3C\equiv0\pmod{17}$ for an element $A+B\theta+C\theta^2$. The congruence defines a lattice of index $17$, so
    \begin{equation}\label{eq:order17}
        R=\Z[t]=\{A+B\theta+C\theta^2\in\OK:B-3C\equiv0\pmod{17}\}.
    \end{equation}

    To pass from the description of $R$ to its unit group, factor the defining polynomial modulo $17$:
    \[
        X^3-3X-1=(X-3)(X-4)(X-10).
    \]
    For $z=A+B\theta+C\theta^2$, equality of its evaluations at $4$ and $10$ is equivalent to $B-3C\equiv0\pmod{17}$. These evaluations are ring homomorphisms $\OK\to\mathbb F_{17}$. On units they are nonzero, and therefore
    \[
        \rho:\OK^\times\longrightarrow\mathbb F_{17}^\times,
        \qquad \rho(z)=z(4)/z(10)
    \]
    is a group homomorphism. Equation~\eqref{eq:order17} shows that its kernel consists of the units of $\OK$ that lie in $R$. Each such unit has its inverse in $R$ by its monic characteristic polynomial, so $\ker\rho=R^\times$.

    Now $\rho(\theta)=4/10=14$ in $\mathbb F_{17}$, with $14^2=9$, $14^4=13$, and $14^8=-1$. Thus $\rho(\theta)$ has order $16$, so $\rho$ is surjective onto $\mathbb F_{17}^\times$ and $[\OK^\times:R^\times]=16$. The exponent columns of $u,v$ in the basis $\theta,\theta+1$ of $\OK^\times/\{\pm1\}$ give
    \[
        [\OK^\times:H]=\left|\det\begin{pmatrix}6&4\\-5&-6\end{pmatrix}\right|=16.
    \]
    Since $H\subseteq R^\times$, equality of the two ambient indices proves $R^\times=H$.
\end{proof}

For this row the two unit indices are $M=16$ and $I_R=1$, while the additive index is $J=17$. A unit of $\OK$ outside $H$ need not belong to $R$. The residue argument proves generation in $R$ without enumerating the sixteen cosets.

\section{A trace condition forcing Galois invariance}\label{sec:trace}
We now give a condition on either unit that makes the joint order Galois invariant. Multiplying a unit in a cubic field by its norm gives a unit of norm one. We consider such normalized units $\psi$ satisfying
\begin{equation}\label{eq:gras-trace}
    \Nm(\psi)=1,\qquad \Tr(\psi+\psi^{-1})=-3.
\end{equation}
For a norm-one cubic unit, $\Nm(1+\psi)=2+\Tr(\psi+\psi^{-1})$, so the trace identity is equivalent to $\Nm(1+\psi)=-1$. The trace identity is preserved by conjugation and inversion. In the next proposition, the field need not be assumed cyclic and the coefficients need not be positive.

\begin{proposition}\label{prop:gras-invariance}
    Let $K$ be a cubic number field, let $u,v\in\OK^\times\setminus\Q$, and let $a,s$ be nonzero integers with $au+sv=1$. Put $n=\Nm(u)$ and $e=\Nm(v)$. If either $nu$ or $ev$ satisfies~\eqref{eq:gras-trace}, then $K/\Q$ is cyclic and $\Z[u,v]$ is Galois invariant.
\end{proposition}
\begin{proof}
    Suppose first that $\psi=nu$ satisfies the condition and put $T_\psi=\Tr(\psi)=nT$. Its minimal polynomial is
    \[
        f_\psi(Z)=Z^3-T_\psi Z^2-(T_\psi+3)Z-1.
    \]
    The transformation $h(z)=-(z+1)/z$ satisfies $z^3f_\psi(h(z))=f_\psi(z)$ and has order three. Neither $0$ nor $-1$ is a root, and $h(\psi)=\psi$ would give a quadratic equation for the degree-three element $\psi$. Thus $\psi,h(\psi),h^2(\psi)$ are the three distinct roots in $K$, proving cyclicity. Choose $\sigma$ with $\sigma(\psi)=h(\psi)$, so $\sigma(u)=-n-u^{-1}$.

    Modulo $s\OK$, the weighted equation gives $u^{-1}\equiv a$, while its conjugate gives $a\sigma(u)\equiv1$. Hence $a^2+na+1\in s\OK\cap\Z=s\Z$, and
    \begin{equation}\label{eq:gras-action}
        q=\frac{a^2+na+1}{s}\in\Z,\qquad
        \sigma(u)=-n-u^{-1},\qquad
        \sigma(v)=q+a\frac vu.
    \end{equation}
    The last equality follows by substituting $v=(1-au)/s$ in $(1-a\sigma(u))/s$. An integral unit's inverse belongs to the order it generates, so both conjugates lie in $R=\Z[u,v]$. Since $\sigma(R)\subseteq R$ and $\sigma^3=1$, all inclusions in $R=\sigma^3(R)\subseteq\sigma^2(R)\subseteq\sigma(R)\subseteq R$ are equalities.

    If $ev$ satisfies the condition, interchange the coefficients and units. For a suitable generator $\tau$ this gives
    \[
        p=\frac{s^2+es+1}{a}\in\Z,\qquad
        \tau(v)=-e-v^{-1},\qquad \tau(u)=p+s\frac uv,
    \]
    and proves the same conclusion.
\end{proof}

\subsection{The trace condition and individual equation orders}
The weighted equation also determines the trace of the normalized unit. If $\psi=nu$ satisfies~\eqref{eq:gras-trace} and $T_\psi=\Tr(\psi)=nT$, then the minimal polynomial of $u$ is
\[
    f_u(X)=X^3-nT_\psi X^2-(T_\psi+3)X-n.
\]
Thus $\Nm(1-au)=a^3f_u(1/a)=1-naT_\psi-a^2(T_\psi+3)-na^3$. Equating this with $\Nm(sv)=es^3$ gives
\begin{equation}\label{eq:trace-parameter}
    a(a+n)T_\psi=1-3a^2-na^3-es^3.
\end{equation}
Under our positive-weight convention, the weights and norm signs therefore determine $T_\psi$ uniquely. The divisibility in~\eqref{eq:gras-action} also gives $s\mid a^2+na+1$, and hence $s\le a^2+na+1$. These conditions are necessary, but the reconstructed second unit must still be checked for integrality.

The same trace parameter determines the index of the order generated by $\psi$. Combining its polynomial discriminant with Proposition~\ref{prop:orderindices} gives the following formula.
\begin{samepage}
    \begin{corollary}\label{cor:trace-index}
        Suppose a pair satisfies~\eqref{eq:weighted}, and either $\psi=nu$ or $\psi=ev$ satisfies~\eqref{eq:gras-trace}. Put $T_\psi=\Tr(\psi)$. Then
        \begin{equation}\label{eq:trace-index}
            \Tr\!\left(\frac{\psi^2-\psi^{-1}}{\mathfrak f_K}\right)
            =\frac{T_\psi^2+3T_\psi+9}{\mathfrak f_K}
            =[\OK:\Z[\psi]]
            =\begin{cases}s^3J,&\psi=nu,\\a^3J,&\psi=ev.\end{cases}
        \end{equation}
        The other unit in the pair does not belong to $\Z[\psi]$.
    \end{corollary}
\end{samepage}
\begin{proof}
    The trace condition gives $f_\psi(Z)=Z^3-T_\psi Z^2-(T_\psi+3)Z-1$, so $\Tr(\psi^{-1})=-T_\psi-3$ and $\Tr(\psi^2)=T_\psi^2+2T_\psi+6$. We also have $\disc(f_\psi)=(T_\psi^2+3T_\psi+9)^2$, while $D_K=\mathfrak f_K^2$ by the conductor--discriminant formula. Since $T_\psi^2+3T_\psi+9>0$, the discriminant-index identity gives the first three equalities. The last equality follows from Proposition~\ref{prop:orderindices}, because multiplying a unit by its norm does not change its equation order.

    If the other unit belonged to $\Z[\psi]$, then $R$ would equal $\Z[\psi]$. This contradicts $[R:\Z[\psi]]>1$, since that index is $s^3$ or $a^3$.
\end{proof}

When $R=\OK$, the opposite weight is therefore the cube parameter in Gras's monogeneity criterion~\cite[\S I]{gras1976}. For a proper order, the extra factor $J$ remains in~\eqref{eq:trace-index}. Galois invariance alone does not imply $J=1$.

\subsection{Adjoining conjugates at prime unit index}\label{sec:prime-unit-index}
The trace condition also gives a criterion for generation after adjoining the conjugates of $u$ and $v$. Let $\sigma$ generate $\operatorname{Gal}(K/\Q)$ and define
\[
    H_{\mathrm{Gal}}=\gp{-1,\sigma^i(u),\sigma^i(v):0\le i\le2},
\]
which is independent of the choice of generator.
\begin{corollary}\label{prop:prime-trace-generation}
    Under the hypotheses of Corollary~\ref{cor:trace-index}, suppose $[\OK^\times:\Z[\psi]^\times]$ is a prime. Then $H_{\mathrm{Gal}}=\OK^\times$.
\end{corollary}
\begin{proof}
    Inversion replaces the trace parameter $T_\psi$ by $-T_\psi-3$ and preserves $\Z[\psi]$. Since $\max\{T_\psi,-T_\psi-3\}\ge-1$, Thomas's Theorem~(3.10)~\cite{thomas1979} applies after this replacement if necessary. Choose the generator $\sigma$ so that $\sigma(\psi)=-(\psi+1)/\psi$. Then
    \[
        \Z[\psi]^\times=\gp{-1,\psi,\sigma(\psi)}
        \subseteq H_{\mathrm{Gal}}.
    \]
    By Corollary~\ref{cor:trace-index}, the other unit lies outside $\Z[\psi]$, and hence outside $\Z[\psi]^\times$. A subgroup containing a subgroup of prime index and an element outside it is the whole group. This proves the assertion.
\end{proof}
The corollary assumes the prime unit index; it does not compute that index from the additive index in~\eqref{eq:trace-index}. It concerns the group generated by all conjugates of the pair and leaves the index of $H$ undetermined.

\section{Unit generation by Minemura's criterion}\label{sec:minemura}
We return to the positive weights $2\le a<s$ and the ratio polynomial~\eqref{eq:P}. Minemura's criterion gives sufficient conditions for $I_R=1$. We express these conditions in terms of $E,F$, then use the coefficient line to locate every possible failure.

Minemura's Theorem~2~\cite{minemura} concerns the norm-one units of the equation order $\Z[\theta]$. The notation $E_\theta^+$ in his theorem denotes norm-one units, not a maximal-order unit group or the totally positive subgroup. We apply the theorem with
\[
    (A,B,C)=(ea,es,c),
\]
and retain the polynomial coefficients $(E,F)$ from \eqref{eq:P}. Minemura's norm identity becomes $e$ times \eqref{eq:normline}, so its right-hand side is one. The generators $e(at+s)$ and $ct$ have norm one. The hypothesis on the sizes of $A$ and $B$ follows from $2\le a<s$, even when $e=-1$.

Writing $\xi=s/a>1$, the remaining sufficient hypotheses are exactly
\begin{align}\label{eq:mine}
    F+E              & \ne1,\qquad F-E\ne1,\notag \\
    |3\xi^2-2E\xi+F| & >|E-3\xi|+1+2a,\notag      \\
    |E-\xi|          & >4a,\qquad |E-2\xi|>5/2.
\end{align}
If these conditions hold, the signed affine units generate the norm-one units of $R=\Z[t]$. Since the degree is odd, adjoining $-1$ gives every unit, hence $R^\times=H$. The order $R$ may be maximal or proper, and need not be Galois invariant. When one of these conditions fails, the theorem gives no conclusion about generation.

For fixed weights and norm signs, the norm identity determines $F$ from $E$. We can therefore locate the coefficients not covered by the criterion on a single line.
\begin{samepage}
    \begin{lemma}\label{lem:intervals}
        For weights and norm signs satisfying \eqref{eq:filters} with $2\le a<s$, put
        \[
            \xi=\frac sa,\qquad q_1=\frac{e-ca^3}{a^2s},\qquad d_1=2a+1,
        \]
        and define
        \[
            L_1=\frac{q_1-\xi-d_1}{\xi-1},\qquad
            U_1=\max\left\{\frac{d_1+q_1+\xi}{\xi+1},\frac{d_1+q_1-\xi}{\xi-1}\right\}.
        \]
        On the coefficient line \eqref{eq:normline}, at least one condition in \eqref{eq:mine} fails precisely when the real coefficient $E$ lies in the union of the intervals
        \begin{align}\label{eq:intervals}
            I_1 & =[2\xi+L_1,2\xi+U_1],\notag \\
            I_2 & =[\xi-4a,\xi+4a],\notag     \\
            I_3 & =[2\xi-5/2,2\xi+5/2]
        \end{align}
        together with the two points
        \begin{equation}\label{eq:points}
            E=\frac{\xi^2+1-q_1}{\xi+1},\qquad E=\frac{\xi^2+1-q_1}{\xi-1}.
        \end{equation}
        The integral coefficients not covered by the criterion are obtained by intersecting this set with the progression
        \begin{equation}\label{eq:residue}
            E\equiv \ell s^{-1}\pmod a,\qquad F=(sE-\ell)/a,
        \end{equation}
        where $\ell$ is defined in \eqref{eq:t}. Nonintegral isolated points contribute nothing.
    \end{lemma}
\end{samepage}
\begin{proof}
    The norm identity gives $F=\xi E-\xi^2+q_1$. Set $z=E-2\xi$. Then the derivative inequality fails precisely when
    \[
        |\xi z-q_1|-|z-\xi|\le d_1.
    \]
    Since $s\ge a+1$,
    \[
        |q_1|\le\frac{a^3+1}{a^2(a+1)}<1,\qquad q_1/\xi<\xi.
    \]
    The left-hand side is the piecewise linear function
    \[
        f(z)=\begin{cases}
            q_1-\xi-(\xi-1)z, & z\le q_1/\xi,       \\
            (\xi+1)z-q_1-\xi, & q_1/\xi\le z\le\xi, \\
            (\xi-1)z+\xi-q_1, & z\ge\xi.
        \end{cases}
    \]
    The function $f$ decreases to $f(q_1/\xi)=q_1/\xi-\xi<0$ and then increases. The left solution of $f(z)=d_1$ is $L_1$. The right solution is $(d_1+q_1+\xi)/(\xi+1)$ if $d_1+q_1\le\xi^2$, and $(d_1+q_1-\xi)/(\xi-1)$ otherwise. Subtracting the first expression from the second gives a quantity with the sign of $d_1+q_1-\xi^2$, so their maximum gives $U_1$ in both cases. This yields $I_1$. The two other strict inequalities give $I_2,I_3$, and $F\pm E=1$ gives \eqref{eq:points}. The endpoints are included because the inequalities in Minemura's theorem are strict. Equation~\eqref{eq:residue} is the integral coefficient line proved earlier.
\end{proof}

The isolated points are redundant. Indeed $|q_1|<1$, $\xi-1\ge1/a$, and
\[
    F+E=1\ \Longrightarrow\ E-\xi=-1+\frac{2-q_1}{\xi+1}\in(-1,1/2),
\]
\[
    F-E=1\ \Longrightarrow\ E-\xi=1+\frac{2-q_1}{\xi-1}\in(1,1+3a).
\]
Both points lie in $I_2$. Thus the coefficients not covered by Minemura's theorem are precisely those in the three intervals that also lie in the coefficient progression.

\begin{theorem}\label{thm:reduction}
    Fix an integer $a\ge2$. To prove $I_R=1$ for every pair satisfying~\eqref{eq:weighted} with $s>a$ and $R\subsetneq\OK$, it suffices to consider the finitely many cubic polynomials satisfying \eqref{eq:filters}, \eqref{eq:intervals}, and \eqref{eq:residue}. Keep only irreducible cubics with positive square discriminant and proper equation order, and prove $\Z[t]^\times=\gp{-1,t,at+s}$ for each of these. Every other coefficient is already covered by Minemura's theorem.
\end{theorem}
\begin{proof}
    There are finitely many choices of weights and norm signs by \eqref{eq:filters}. For each choice, the intervals are bounded, so their intersection with the arithmetic progression is finite. Lemma~\ref{lem:ratio} gives the correspondence with weighted units, and Minemura's theorem proves generation for all remaining coefficients.
\end{proof}
This procedure does not require solving the entire square-discriminant problem for any weight. For example, $I_2$ has length $8a$, so it contains at most nine entries of the progression modulo $a$; $I_3$ has length five. Enumerating each interval separately avoids scanning the gaps between them when $s/a$ is large. The same formulas apply to every $a$, but the full set of parameters is unbounded when $a$ varies.

\section{Weights one away from a cube}\label{sec:extreme}
When the second weight is one away from the cube of the first, the intervals in Section~\ref{sec:minemura} leave only the polynomial families described below. The proposition does not decide which of these polynomials have square discriminant.

\begin{proposition}\label{prop:extreme-reduction}
    Suppose~\eqref{eq:weighted} holds with $a\ge3$ and $s=a^3\pm1$. Then $s=a^3+e$, $n=-e$, and the ratio polynomial has the form
    \begin{equation}\label{eq:extremeP}
        P_{a,j,e}(X)=X^3+(a^2+aj)X^2+\bigl((a^3+e)j-ea\bigr)X+1,
        \qquad j\in\Z.
    \end{equation}
    Put
    \[
        \mathcal J_{1}=\{-3,-2,-1,0,1,2,3,4\},\qquad
        \mathcal J_{-1}=\{-4,-3,-2,-1,0,1,2,3\}.
    \]
    If $R^\times\ne H$, then $j\in\mathcal J_e\setminus\{0\}$ and $j\ne a$. Thus any failure of generation for these weights must occur for one of fourteen choices of $(e,j)$, with $a$ still variable.
\end{proposition}
\begin{proof}
    Write $s=a^3+\varepsilon$, where $\varepsilon\in\{\pm1\}$. The norm sign $e=\Nm(v)$ satisfies $a\mid es^3-1$, so $a\mid e\varepsilon-1$. Since $a\ge3$, this forces $e=\varepsilon$. Similarly, $s\mid na^3-1$ forces $n=-e$, because $a^3\equiv-e\pmod s$ and $s>2$. Thus $s=a^3+e$ and $c=ne=-1$.

    The coefficient line now has $\ell=(s^2-1)/a=a^2(s+e)$ in~\eqref{eq:t}. Thus $E\equiv0\pmod a$, and writing $E=a^2+aj$ gives $F=sj-ea$, proving~\eqref{eq:extremeP}. Conversely, these coefficients satisfy the norm identity; irreducibility and positive square discriminant are the remaining conditions for a cyclic pair.

    Set $\xi=s/a=a^2+e/a$, $q_1=1/a^2$, and $z=E-2\xi=a(j-a)-2e/a$. Failure of Minemura's derivative inequality implies
    \[
        |\xi z-q_1|\le|z-\xi|+2a+1,
        \qquad
        |z|\le\frac{2a+1+\xi+q_1}{\xi-1}
        <1+\frac{2a+3}{a^2-1-1/a}.
    \]
    If $j\ne a$, then $|z|\ge a-2/a$. For $a\ge3$,
    \[
        a-2/a>1+\frac{2a+3}{a^2-1-1/a}.
    \]
    Indeed, after multiplication by the positive denominator and by $a^2$, the difference has numerator
    \[
        G(a)=a^5-a^4-5a^3-3a^2+3a+2.
    \]
    Substitution $a=b+3$ gives
    \[
        G(b+3)=b^5+14b^4+73b^3+168b^2+147b+11>0
        \qquad(b\ge0).
    \]
    Therefore derivative failure is possible only at $j=a$.

    The inequality $|E-\xi|>4a$ fails exactly when $|j-e/a^2|\le4$, whose integral solutions are $\mathcal J_e$. For $a\ge4$, failure of $|E-2\xi|>5/2$ again forces $j=a$, since $a-2/a>5/2$. At $a=3$ the additional possible value is $j=4$ for $e=1$, or $j=2$ for $e=-1$; both already lie in $\mathcal J_e$. The equalities $F\pm E=1$ lie in the failed interval $|E-\xi|\le4a$ by Section~\ref{sec:minemura}. Thus all failures of Minemura's hypotheses lie in $\mathcal J_e\cup\{a\}$.

    Two parameters can be removed without any finite search. At $j=a$ the discriminant is
    \[
        4a^6-27,
        \qquad (2a^3-1)^2<4a^6-27<(2a^3)^2
        \quad(a\ge3),
    \]
    and at $j=0$ it is
    \[
        -3a^6-14ea^3-27=-\frac{(3a^3+7e)^2+32}{3}<0.
    \]
    Neither can define a cyclic cubic field. Every other parameter outside $\mathcal J_e\setminus\{0\}$ satisfies Minemura's sufficient criterion, proving the assertion.
\end{proof}

The fourteen choices still allow $a$ to vary, so the proposition does not prove finiteness of the cyclic pairs in these families. Failure of Minemura's inequalities leaves generation undecided. A further argument is needed either to exclude the remaining cyclic pairs or to prove $R^\times=\gp{-1,u,v}$ for each pair.

\section{Two open questions}\label{sec:global}
The first question asks whether $I_R=1$ for every proper joint order arising from the weighted equation.
\begin{conjecture}\label{conj:generation}
    For every pair satisfying~\eqref{eq:weighted}, if $R\subsetneq\OK$ then $R^\times=H$.
\end{conjecture}
Note that Proposition~\ref{prop:consecutive-generation} proves the conjecture for all consecutive pairs. Minemura's theorem settles every pair whose coefficient $E$ lies outside the intervals in Section~\ref{sec:minemura}, and Theorem~\ref{thm:reduction} leaves finitely many cases for each fixed first weight. Proposition~\ref{prop:extreme-reduction} reduces possible failures at $s=a^3\pm1$ to fourteen polynomial families.

The second question asks for a bound on the weights that holds over all cyclic cubic fields.
\begin{conjecture}\label{conj:weights}
    There is an absolute constant $C$ such that every pair satisfying~\eqref{eq:weighted}, over every cyclic cubic field, has $s\le C$.
\end{conjecture}
This conjecture is a uniform weighted analogue of the fixed-field finiteness results for $N_K=\{m\in\Z_{>0}:m=\varepsilon+\delta,\ \varepsilon,\delta\in\OK^\times\}$ of Tinkov\'a--Visser--Yatsyna~\cite[Theorem~1.1]{tinkovavisseryatsyna} and Visser--Ziegler~\cite[Corollary~1]{VisserZiegler2026}, with the bound required to hold as $K$ varies over all cyclic cubic fields.

Together with Corollary~\ref{cor:fixed}, Conjecture~\ref{conj:weights} would imply finiteness of all solutions of~\eqref{eq:weighted} up to simultaneous $\Q$-isomorphism.

\bibliographystyle{alpha}
\begingroup
\raggedright
\bibliography{references}

@article{evertsegyory1988,
  author  = {Evertse, Jan-Hendrik and Gy{\H{o}}ry, K{\'a}lm{\'a}n},
  title   = {On the numbers of solutions of weighted unit equations},
  journal = {Compositio Mathematica},
  volume  = {66},
  number  = {3},
  year    = {1988},
  url     = {https://www.numdam.org/item/CM_1988__66_3_329_0/}
}

@article{vukusicziegler2022,
  author  = {Vukusic, Ingrid and Ziegler, Volker},
  title   = {On a family of unit equations over simplest cubic fields},
  journal = {Journal de Th{\'e}orie des Nombres de Bordeaux},
  volume  = {34},
  number  = {3},
  year    = {2022},
  doi     = {10.5802/jtnb.1223}
}

@article{brindzagyory,
  author  = {Brindza, B{\'e}la and Gy{\H{o}}ry, K{\'a}lm{\'a}n},
  title   = {On unit equations with rational coefficients},
  journal = {Acta Arithmetica},
  volume  = {53},
  year    = {1990},
  pages   = {367--388},
  doi     = {10.4064/aa-53-4-367-388},
  note    = {\url{https://matwbn.icm.edu.pl/ksiazki/aa/aa53/aa5345.pdf}}
}

@article{gras1976,
  author  = {Gras, Marie-Nicole},
  title   = {Lien entre le groupe des unit{\'e}s et la monog{\'e}n{\'e}it{\'e} des corps cubiques cycliques},
  journal = {Publications math{\'e}matiques de la Facult{\'e} des sciences de Besan{\c{c}}on. Th{\'e}orie des nombres},
  year    = {1976},
  pages   = {1--20},
  doi     = {10.5802/pmb.a-9},
  note    = {No.~1, article no.~4; \url{https://www.numdam.org/articles/10.5802/pmb.a-9/}}
}

@article{komatsu,
  author  = {Komatsu, Toru},
  title   = {Pairs of cyclic cubic units with rational difference},
  journal = {International Journal of Number Theory},
  volume  = {21},
  number  = {1},
  year    = {2025},
  pages   = {35--52},
  doi     = {10.1142/S1793042125500034},
  note    = {\url{https://doi.org/10.1142/S1793042125500034}}
}

@article{levin,
  author  = {Levin, Aaron},
  title   = {Variations on a theme of {Runge}: effective determination of integral points on certain varieties},
  journal = {Journal de Th{\'e}orie des Nombres de Bordeaux},
  volume  = {20},
  year    = {2008},
  pages   = {385--417},
  doi     = {10.5802/jtnb.634},
  note    = {\url{https://www.numdam.org/item/JTNB_2008__20_2_385_0/}}
}

@article{minemura,
  author  = {Minemura, Kenji},
  title   = {On Totally Real Cubic Orders Whose Unit Groups are of Type {$\langle a\theta+b,c\theta+d\rangle$}},
  journal = {Tokyo Journal of Mathematics},
  volume  = {27},
  year    = {2004},
  pages   = {393--410},
  doi     = {10.3836/tjm/1244208398},
  note    = {\url{https://doi.org/10.3836/tjm/1244208398}}
}

@article{thomas1979,
  author  = {Thomas, Emery},
  title   = {Fundamental units for orders in certain cubic number fields},
  journal = {Journal f{\"u}r die reine und angewandte Mathematik},
  volume  = {310},
  year    = {1979},
  pages   = {33--55},
  doi     = {10.1515/crll.1979.310.33},
  note    = {\url{https://doi.org/10.1515/crll.1979.310.33}}
}

@article{tinkovavisseryatsyna,
  author  = {Tinkov{\'a}, Magdal{\'e}na and Visser, Robin and Yatsyna, Pavlo},
  title   = {Sums of two units in number fields},
  journal = {Mathematische Zeitschrift},
  volume  = {312},
  year    = {2026},
  pages   = {Paper No.~39},
  doi     = {10.1007/s00209-025-03919-z},
  note    = {\url{https://doi.org/10.1007/s00209-025-03919-z}}
}

@article{thomasvasquez,
  author  = {Thomas, Emery and Vasquez, Alphonse T.},
  title   = {Diophantine equations arising from cubic number fields},
  journal = {Journal of Number Theory},
  volume  = {13},
  year    = {1981},
  pages   = {398--414},
  doi     = {10.1016/0022-314X(81)90023-8},
  note    = {\url{https://doi.org/10.1016/0022-314X(81)90023-8}}
}

@misc{VisserZiegler2026,
  author = {Robin Visser and Volker Ziegler},
  title  = {On integers that are representable as the sum of two units},
  year   = {2026},
  note   = {Preprint, arXiv:2608.15903. \url{https://arxiv.org/abs/2608.15903}}
}
\endgroup
\end{document}